\documentclass{article}
\usepackage[english]{babel}
\usepackage{latexsym,amssymb,amsmath}
\usepackage{amsfonts,amssymb}
\usepackage{euscript}
\newcounter{theorem} 
\newcounter{lemma} 

\,
\,

\begin{document}

\noindent {\bf \Large Cauchy theorem for a surface integral in
commutative algebras}
\vspace{2mm}\\
\noindent {\bf \large S.\,A.~Plaksa and V.\,S.~Shpakivskyi}\\
\vspace{5mm}


\centerline{\small\sl Dedicated to memory of Professor Promarz
M.~Tamrazov}

\vskip 2mm

\begin{abstract}
We prove an analogue of the Cauchy integral theorem for
hyperholomorphic functions given in three-dimensional domains with
non piece-smooth boundaries and taking values in an arbitrary
finite-dimensional commutative associative Banach algebra.
\end{abstract}

\vspace{2mm}

\section{Introduction}

The Cauchy integral theorem is a fundamental result of the
classical complex analysis in the complex plane $\mathbb{C}$: if
the boundary $\partial D$ of a domain $D\subset\mathbb{C}$ is a
closed Jordan rectifiable curve, and a function $F\colon
\overline{D}\longrightarrow \mathbb{C}$ is continuous in the
closure   $\overline{D}$ of $D$ and is holomorphic in $D$,
then\,\, $$\int\limits_{\partial D}F(z)dz=0\,.$$

Developing hypercomplex analysis in both commutative and
noncommutative algebras needs similar general analogues of the
Cauchy integral theorem for several-dimensional spaces.

It is well known that in the case where a simply connected domain
has a closed piece-smooth boundary, spatial analogues of the
Cauchy integral theorem can be obtained with using the classical
Gauss -- Ostrogradskii formula, if a given function has
specifically continuous partial derivatives of the first order up
to the boundary. In such a way analogues of the Cauchy integral
theorem are proved in the quaternion algebra (see, e.g.,
\cite[p.~66]{Kravchenko-Shapiro}) and in Clifford algebras (see,
e.g., \cite[p.~52]{Brakx-Delan-Som}).

Generalizations of the Cauchy integral theorem have relations to
weakening requirements to the boundary or the given function.
Usually, such generalizations are based on generalized Gauss --
Ostrogradskii -- Green -- Stokes formula (see, e.g.,
\cite{Federer,Harrison-Norton}) under the condition of continuity
of partial derivatives of the given function, but for extended
classes of surfaces of integration; see, e.g.,
\cite{Abreu-Bory-99,Abreu-Bory-Pena}, where rectifiable or regular
surfaces are considered. In the papers \cite{Sadberi,Gerus-2011}
the continuity of partial derivatives is changed by a
differentiability of components of the given function taking
values in the quaternion algebra. Note that the boundary of domain
is remained piece-smooth in \cite{Gerus-2011}.

In this paper we prove an analogue of the Cauchy integral theorem
for functions taking values in an arbitrary finite-dimensional
commutative associative algebra. Similarly to the paper
\cite{Gerus-2011}, we weaken requirements to functions given in a
domain of three-dimensional space. At the same time, the functions
can be given in a domain with non piece-smooth boundary.

\section{Quadrable surfaces}

A set $\Sigma$ is called a \textit{surface} in the real space
$\mathbb{R}^{3}$ if $\Sigma$ is a homeomorphic image of the square
$G:=[0,1]\times [0,1]$ (cf., e.g., \cite[p.~24]{Rado}).

By $\Sigma^\varepsilon$ we denote $\varepsilon$-neighborhood of
the surface $\Sigma$, i.e. the set
$\Sigma^\varepsilon:=\{(x,y,z)\in\mathbb{R}^3:
\sqrt{(x-x_1)^2+(y-y_1)^2+(z-z_1)^2}\leq\varepsilon,
 \, (x_1,y_1,z_1)\in \Sigma\}.$

The \textit{Fr\'{e}chet distance} $d(\Sigma,\Lambda)$ between the
surfaces $\Sigma$ and $\Lambda$ is called the infimum of real
numbers $\varepsilon$, for which the relations $\Sigma\subset
\Lambda^\varepsilon$, $\Lambda\subset \Sigma^\varepsilon$ are
fulfilled  (see, e.g., \cite{Freshe}). A sequence of polyhedral
surfaces $\Lambda_n$ \textit{converges uniformly} to the surface
$\Sigma$, if $d(\Lambda_n,\Sigma)\rightarrow0$ as
$n\rightarrow\infty$ (see, e.g., \cite[p. 121]{Rado}).

The \textit{Lebesgue area} of a surface $\Sigma$ is
$$\mathfrak{L}(\Sigma):=\inf\liminf\limits_{n\rightarrow\infty}
\mathfrak{L}(\Lambda_n),$$ where the infimum is taken for all
sequences $\Lambda_n$ convergent uniformly to $\Sigma$ (see, e.g.,
\cite[p. 468]{Rado}), and $\mathfrak{L}(\Lambda_n)$ is the area of
polyhedral surface $\Lambda_n$.

Let a surface $\Sigma$ have the finite Lebesgue area, i.e.
$\mathfrak{L}(\Sigma)<\infty$. Then by the L.~Cesari theorem
\cite[p.~544]{Cesari}, there exists a surface parameterization
$$
\Sigma=\left\{f(u,v):=\big(x(u,v),\,y(u,v),\,z(u,v)\big):(u,v)\in
G\right\}
$$
such that the Jacobians
\begin{equation}\label{Jacobiany}
A:=\frac{\partial y}{\partial u}\frac{\partial z}{\partial
v}-\frac{\partial y}{\partial v}\frac{\partial z}{\partial u},
\quad B:=\frac{\partial z}{\partial u}\frac{\partial x}{\partial
v}-\frac{\partial z}{\partial v}\frac{\partial x}{\partial
u},\quad  C:=\frac{\partial x}{\partial u}\frac{\partial
y}{\partial v}-\frac{\partial x}{\partial v}\frac{\partial
y}{\partial u}
\end{equation}
exist a.e. in the square $G$ and
\begin{equation}\label{riv-pl-leb}
\mathfrak{L}(\Sigma)=\int\limits_{G}\sqrt{A^2+B^2+C^2}\,dudv
\end{equation}
(here and in what follows, all integrals are understood as
Lebesgue integrals).

In the case where $\mathfrak{L}(\Sigma)<\infty$ and the equality
(\ref{riv-pl-leb}) holds for the given parameterization of
$\Sigma$, we shall say that a surface $\Sigma$ is {\it quadrable}.

Let us formulate certain sufficient conditions for a surface
$\Sigma$ be quadrable.

 \begin{itemize}
\item[1.] If $\Sigma$ is a rectifiable surface, then it follows from
 \cite[IV.4.28, IV.4.1~(e)]{Rado} that $\Sigma$ is
quadrable.

\item[2.] Let the components $x(u,v),\,y(u,v),\,z(u,v)$ of mapping
$f$ be absolutely continuous in the sense of Tonelli (see., e.~g.,
\cite[p.~169]{Saks}). Let, furthermore, in Jacobians $A,B,C$ of
mapping $f$ in every of the products $\frac{\partial y}{\partial
u}\frac{\partial z}{\partial v}$, $\frac{\partial y}{\partial
v}\frac{\partial z}{\partial u}$, $\frac{\partial z}{\partial
u}\frac{\partial x}{\partial v}$, $\frac{\partial z}{\partial
v}\frac{\partial x}{\partial u}$, $\frac{\partial x}{\partial
u}\frac{\partial y}{\partial v}$, $\frac{\partial x}{\partial
v}\frac{\partial y}{\partial u}$, one partial derivative belong to
the class $L_p(G)$ of functions integrable to the $p$th power on
$G$ and the other partial derivative belong to $L_q(G)$, where
$\frac{1}{p}+\frac{1}{q}=1$. Then $\Sigma$ is quadrable (see
\cite[V.2.26]{Rado}). Note that for a rectifiable surface
$\Sigma$, components $x(u,v),\,y(u,v),\,z(u,v)$ of mapping $f$ are
absolutely continuous in the sense of Tonelli (see, e.g.,
\cite[p.~169]{Saks}).

\item[3.] If two components of the mapping $f(u,v)$ are Lipschitz
functions and the third component is absolutely continuous in the
sense of Tonelli, then $\Sigma$ is quadrable (see
\cite[V.2.28]{Rado}).
\end{itemize}

\section{Surface integrals}

In what follows, we understand the {\it closed surface}
$\Gamma\subset\mathbb{R}^3$ as an image of a sphere under
homeomorphic mapping which maps {\it at least one circle onto a
rectifiable curve}. In other words, the closed surface $\Gamma$ is
the union of two surfaces $\Gamma_1$, $\Gamma_2$ for which
$\Gamma_1\cap\Gamma_2=:\gamma$ is a closed Jordan rectifiable
curve. Let the surfaces $\Gamma_1$, $\Gamma_2$ be parametrically
definable:
$$
\Gamma_1=\left\{f_1(u,v):=\big(x_1(u,v),\,y_1(u,v),\,z_1(u,v)\big):(u,v)\in
G\right\},
$$
$$
\Gamma_2=\left\{f_2(u,v):=\big(x_2(u,v),\,y_2(u,v),\,z_2(u,v)\big):(u,v)\in
G\right\}.$$

A closed surface $\Gamma$ is called {\it quadrable}\/ if the
surfaces $\Gamma_1$ and $\Gamma_2$ are quadrable.

For a closed quadrable surface $\Gamma$ and a continuous function
$F:\Gamma\rightarrow\mathbb{R}$, we define integrals on $\Gamma$
by the equalities
$$\int\limits_{\Gamma}F(x,y,z)\,dydz:=\int\limits_{G}F\Big(x_1(u,v),y_1(u,v),z_1(u,v)\Big)A_1
\,dudv-$$
\begin{equation}\label{int-1}
-\int\limits_{G}F\Big(x_2(u,v),y_2(u,v),z_2(u,v)\Big)A_2\,dudv,
\end{equation}
$$\int\limits_{\Gamma}F(x,y,z)\,dzdx:=\int\limits_{G}F\Big(x_1(u,v),y_1(u,v),z_1(u,v)\Big)B_1
\,dudv-$$
\begin{equation}\label{int-2}
-\int\limits_{G}F\Big(x_2(u,v),y_2(u,v),z_2(u,v)\Big)B_2\,dudv,
\end{equation}
$$\int\limits_{\Gamma}F(x,y,z)\,dxdy:=\int\limits_{G}F\Big(x_1(u,v),y_1(u,v),z_1(u,v)\Big)C_1
\,dudv-$$
\begin{equation}\label{int-3}
-\int\limits_{G}F\Big(x_2(u,v),y_2(u,v),z_2(u,v)\Big)C_2\,dudv
\end{equation}
with the Jacobians $A_k,B_k,C_k$ of mapping $f_k$ of the form
(\ref{Jacobiany}) for $k=1,2$.

It is easy to check up that the definitions (\ref{int-1})
--- (\ref{int-3}) are
correct. Indeed, values of integrals on the right-hand sides of
equalities (\ref{int-1}) --- (\ref{int-3}) are the same for all
parameterizations $f_1$, $f_2$ for which the areas
$\mathfrak{L}(\Gamma_1)$, $\mathfrak{L}(\Gamma_2)$ are expressed
by the equalities of the form (\ref{riv-pl-leb}), and values of
integrals on the left-hand sides of equalities (\ref{int-1}) ---
(\ref{int-3}) do not depend on a choice of a rectifiable curve
$\gamma$ which divides $\Gamma$ into two parts.\vskip 2mm

{\bf Lemma 1.} {\it If $\Gamma$ is a closed quadrable surface,
then}
\begin{equation}\label{lem-3-int-rivni-nulyu}
\int\limits_{\Gamma}dydz=\int\limits_{\Gamma}dzdx=\int\limits_{\Gamma}dxdy=0.
\end{equation}

\textbf{\textit{Proof.}} By definition,
\begin{equation}\label{lem-3-1}
\int\limits_{\Gamma}dydz=\int\limits_{G}A_1\,dudv-\int\limits_{G}A_2\,dudv.
\end{equation}
It follows from the Rad\'o results \cite[V.2.64~$(iii)$,
IV.4.21~$(iii_3)$]{Rado} that for the surfaces $\Gamma_1$,
$\Gamma_2$ the following equalities are true:
\begin{equation}\label{lem-3-2}
\int\limits_{G}A_k\,dudv=\int\limits_{\partial G}ydz,\qquad k=1,2,
\end{equation}
where the integral on the right-hand side is understood as a
Lebesgue -- Stieltjes integral and is took along the boundary
$\partial G$ of the square $G$ into a positive direction. Now, we
obtain from the equalities (\ref{lem-3-1}), (\ref{lem-3-2}) that
the first integral of (\ref{lem-3-int-rivni-nulyu}) is equal to
zero. The other equalities (\ref{lem-3-int-rivni-nulyu}) are
proved by analogy. Lemma is proved.

\section{Hyperholomorphic functions in a
commutative Banach algebra}

Let $\mathbb{A}$ be a commutative associative Banach algebra over
the field of complex numbers $\mathbb{C}$ with the basis
$\{e_{k}\}_{k=1}^{n}$, $3\leq n<\infty$.

Let us single out the linear span $E_{3}:=\{\zeta=x+ye_{2}+ze_{3}:
x, y, z\in\mathbb{R}\}$ generated by the vectors $e_1,e_2,e_3$.
Associate with a set $\Omega\subset{\mathbb R}^3$ the set
$\Omega_{\zeta}:=\{\zeta=x+ye_{2}+ze_{3}: (x, y, z)\in \Omega\}$
in $E_{3}$.

Consider a function $\Psi:\Omega_{\zeta}\rightarrow\mathbb{A}$ of
the form
\begin{equation}\label{Fi}
\Psi(\zeta)=\sum\limits_{k=1}^{n}U_{k}(x,y,z)e_{k}+i\sum\limits
_{k=1}^{n}V_{k}(x,y,z)e_{k},
\end{equation}
where $(x,y,z)\in\Omega$ and $U_k : \Omega\rightarrow\mathbb{R}$,
$V_k : \Omega\rightarrow\mathbb{R}$.

We shall say that a function of the form (\ref{Fi}) is
\textit{hyperholomorphic}\/ in a domain $\Omega_\zeta$ if its
real-valued components $U_k, V_k$ are differentiable in $\Omega$
and the following equality is fulfilled in every point of
$\Omega_\zeta$:
\begin{equation}\label{GiperGol}
\frac{\partial\Psi}{\partial x}\,e_1+\frac{\partial\Psi}{\partial
y}\,e_{2}+\frac{\partial\Psi}{\partial z}\,e_{3}=0.
\end{equation}

In the scientific literature the different denominations are used
for functions satisfying equations of the form (\ref{GiperGol}).
For example, in the papers \cite{Sadberi,Colombo,Spr} they are
called regular functions, and in the papers
\cite{Brakx-Delan-Som,Bernstein,Ryan} they are called monogenic
functions. We use the terminology of the papers
\cite{Kravchenko-Shapiro, Spros,Gerus-2011}.

\section{Auxiliary results}

Let $\Omega$ be a bounded closed set in $\mathbb{R}^{3}$. For a
continuous function $\Psi:\Omega_{\zeta}\rightarrow\mathbb{A}$ of
the form (\ref{Fi}), we define a volume integral by the equality
$$\int\limits_{\Omega_{\zeta}}\Psi(\zeta)dxdydz:=\sum\limits_{k=1}^{n}
e_{k}\int\limits_{\Omega}U_{k}(x,y,z)dxdydz+
i\sum\limits_{k=1}^{n}e_{k}\int\limits_{\Omega}V_{k}(x,y,z)dxdydz.$$

Let $\Gamma$ be a closed quadrable surface in $\mathbb{R}^{3}$.
For a continuous function $\Psi:\Gamma_{\zeta}\rightarrow
\mathbb{A}$ of the form (\ref{Fi}), where $(x,y,z)\in\Gamma$ and
$U_k : \Gamma\rightarrow\mathbb{R}$, $V_k :
\Gamma\rightarrow\mathbb{R}$, we define a surface integral on
$\Gamma_{\zeta}$ with the differential form
$\sigma:=dydze_1+dzdxe_{2}+dxdye_{3}$ by the equality
 $$\int\limits_{\Gamma_{\zeta}}\Psi(\zeta)\sigma:=
 \sum\limits_{k=1}^{n}e_1e_{k}\int\limits_{\Gamma}U_{k}(x,y,z)dydz+
\sum\limits_{k=1}^{n}e_{2}e_{k}\int\limits_{\Gamma}U_{k}(x,y,z)dzdx+$$\vspace{1mm}
$$+\sum\limits_{k=1}^{n}e_{3}e_{k}\int\limits_{\Gamma}U_{k}(x,y,z)dxdy
+i\sum\limits_{k=1}^{n}e_1e_{k}\int\limits_{\Gamma}V_{k}(x,y,z)dydz+$$\vspace{1mm}
$$+i\sum\limits_{k=1}^{n}e_{2}e_{k}\int\limits_{\Gamma}V_{k}(x,y,z)dzdx+
i\sum\limits_{k=1}^{n}e_{3}e_{k}\int\limits_{\Gamma}V_{k}(x,y,z)dxdy,$$\vspace{2mm}
where the integrals on the right-hand side of equality are defined
by the equalities (\ref{int-1}) --- (\ref{int-3}).

The next lemma is a result of Lemma 1 and the definition of
$\sigma$.\vskip 2mm

{\bf Lemma 2.} \textit{If $\Gamma$ is a closed quadrable surface,
then}
\begin{equation}\label{int-sigma}
\int\limits_{\Gamma_\zeta}\sigma=0.
\end{equation}

Let us introduce the Euclidian norm
$\|a\|:=\bigl(\sum\limits_{k=1}^n|a_k|^2\bigr)^{1/2}$ in the
algebra $\mathbb{A}$, where $a=\sum\limits_{k=1}^na_ke_k$ and
$a_k\in\mathbb{C}$ for $k=\overline{1,n}$.

Let $\Gamma$ be a closed quadrable surface in $\mathbb{R}^{3}$.
For a continuous function $U : \Gamma_{\zeta}\rightarrow
\mathbb{R}$, we define a surface integral on $\Gamma_{\zeta}$ with
the differential form $\|\sigma\|$ by the equality \vspace{2mm}
$$\int\limits_{\Gamma_\zeta}U(xe_1+ye_{2}+ze_{3})\|\sigma\|:=$$
$$=\int\limits_{G}U\Big(x_1(u,v)e_1
+y_1(u,v)e_{2}+z_1(u,v)e_{3}\Big)
\sqrt{A_1^2+B_1^2+C_1^2}\,dudv-$$
$$-\int\limits_{G}U\Big(x_2(u,v)e_1
+y_2(u,v)e_{2}+z_2(u,v)e_{3}\Big) \sqrt{A_2^2+B_2^2+C_2^2}\,dudv.
$$\vskip 2mm

{\bf Lemma 3.} \textit{If $\Gamma$ is a closed quadrable surface
and a function $\Psi:\Gamma_{\zeta}\rightarrow \mathbb{A}$ is
continuous, then
$$\Biggr\|\int\limits_{\Gamma_{\zeta}}\Psi(\zeta)\sigma\Biggr\|\leq
 3nM \int\limits_{\Gamma_{\zeta}}\|\Psi(\zeta)\|\|\sigma\|$$
with $M:=\max\limits_{1\leq m,s\leq n}\|e_me_s\|$. }\vskip 2mm

\enlargethispage{2\baselineskip}

\textbf{\textit{Proof.}} Using the representation (\ref{Fi}),
where $(x,y,z)\in\Gamma$, we obtain
 $$\Biggl\|\int\limits_{\Gamma_{\zeta}}\Psi(\zeta)\sigma\Biggr\|\le
 \sum\limits_{k=1}^{n}\|e_{1}e_{k}\|\int\limits_{\Gamma}\bigl|U_{k}(x,y,z)+iV_{k}(x,y,z)
 \bigr|\,dydz+$$\vspace{1mm}
$$+\sum\limits_{k=1}^{n}\|e_{2}e_{k}\|\int\limits_{\Gamma}\bigl|U_{k}(x,y,z)
+iV_{k}(x,y,z)\bigr|\,dzdx+$$\vspace{1mm}
$$+\sum\limits_{k=1}^{n}\|e_{3}e_{k}\|\int\limits_{\Gamma}\bigl|U_{k}(x,y,z)
+iV_{k}(x,y,z)\bigr|\,dxdy\le 3nM
\int\limits_{\Gamma_{\zeta}}\|\Psi(\zeta)\|\|\sigma\|\,.
$$
Lemma is proved. \vskip 2mm

If a simply connected domain $\Omega\subset\mathbb{R}^{3}$ have a
closed piece-smooth boundary $\partial \Omega$ and a function
$\Psi:\Omega_{\zeta}\rightarrow\mathbb{A}$ is continuous together
with partial derivatives of the first order up to the boundary
$\partial \Omega_{\zeta}$, then the following equality follows
from the classical Gauss -- Ostrogradskii formula:
\begin{equation}\label{form-Ostrogradsky}\int\limits_{\partial
\Omega_{\zeta}}\Psi(\zeta)\sigma=\int\limits_{\Omega_{\zeta}}\left(\frac{\partial
\Psi}{\partial x}\,e_1+\frac{\partial \Psi}{\partial
y}\,e_{2}+\frac{\partial \Psi}{\partial z}\,e_{3}\right)dxdydz\,.
\end{equation}

We prove the next theorem similarly to the proof of Theorem 9
\cite{Sadberi} and Theorem 1 \cite{Gerus-2011}, where functions
taking values in the quaternion algebra was considered.\vskip 2mm

{\bf Theorem 1.} \textit{Let $\partial P$ be the boundary of a
closed cube $P$ that is contained in a domain $\Omega$ and a
function $\Psi:\Omega_{\zeta}\rightarrow\mathbb{A}$ be
hyperholomorphic in the domain $\Omega_{\zeta}$. Then the
following equality holds:}
$$\int\limits_{\partial P_{\zeta}}\Psi(\zeta)\sigma=0.$$

\textbf{\textit{Proof.}}\,\, Suppose that\,\,
$\bigl\|\int_{\partial
P_{\zeta}}\Psi(\zeta)\sigma\bigr\|=K.$ 
\vskip 1mm

Denote by $S$ the area of surface $\partial P$. Divide $P$ into
$8$ equal cubes and denote by $P^1$ such a cube, for which \,\,
$\bigl\|\int_{\partial P^1_{\zeta}}\Psi(\zeta)\sigma\bigr\|\geq
K/8$. Clearly, the surface $\partial P^1$ have the area
$S/4$.\vskip 1mm

Continuing this process, we obtain a sequence of embedded cubes
$P^m$ with the areas $S/4^m$ of the surfaces $\partial P^m$, that
satisfies the inequalities
\begin{equation}\label{ocinka-1}
\Biggr\|\int\limits_{\partial
P^m_{\zeta}}\Psi(\zeta)\sigma\Biggr\|\ge K/8^m.
\end{equation}

By the Cantor principle, there exists the unique point
$\zeta_0:=x_0e_1+y_0e_2+z_0e_3$ common for all cubes $P^m$.
Inasmuch as the function $\Psi$ is of the form (\ref{Fi}) and the
real-valued components $U_k, V_k$ are differentiable in $\Omega$,
in a neighbourhood of the point $\zeta_0$ we have the expansion
$$\Psi(\zeta)=\Psi(\zeta_0)+\Delta x\frac{\partial
\Psi(\zeta_0)}{\partial x}+\Delta y\frac{\partial
\Psi(\zeta_0)}{\partial y}+\Delta z\frac{\partial
\Psi(\zeta_0)}{\partial z}+\delta(\zeta,\zeta_0)\rho,$$ where
$\Delta x:=x-x_0$, $\Delta y:=y-y_0$, $\Delta z:=z-z_0$, and
$\delta(\zeta,\zeta_0)$ is an infinitesimal function as
$\rho:=\|\zeta-\zeta_0\|\rightarrow0$.

Therefore, for all sufficiently small cubes, we have
$$\int\limits_{
\partial P^m_{\zeta}}\Psi(\zeta)\sigma=\Psi(\zeta_0)\int\limits_{
\partial P^m_{\zeta}}\sigma+\frac{\partial \Psi(\zeta_0)}{\partial
x}\int\limits_{\partial P^m_{\zeta}}\Delta x\sigma+\frac{\partial
\Psi(\zeta_0)}{\partial y}\int\limits_{ P^m_{\zeta}}\Delta
y\sigma+$$
$$\frac{\partial \Psi(\zeta_0)}{\partial z}\int\limits_{
\partial P^m_{\zeta}}\Delta z\sigma+\int\limits_{
\partial P^m_{\zeta}}\delta(\zeta,\zeta_0)\rho\,\sigma=\sum_{r=1}^5I_r.$$

By the formula (\ref{form-Ostrogradsky}),\,\,  $I_1=0$. Using
(\ref{form-Ostrogradsky}) and taking into account the equality
(\ref{GiperGol}), we obtain
$$I_2+I_3+I_4=\frac{\partial \Psi(\zeta_0)}{\partial
x}\,e_1V_m+\frac{\partial \Psi(\zeta_0)}{\partial
y}\,e_2V_m+\frac{\partial \Psi(\zeta_0)}{\partial z}\,e_3V_m=0,$$
where by $V_m$ we have denoted the volume of cube $P^m$.

Note that for an arbitrary $\varepsilon>0$ there exists the number
$m_0$ such that the inequality
$\|\delta(\zeta,\zeta_0)\|<\varepsilon$ is fulfilled for all cubes
$P^m$ with $m>m_0$. Note also that $\rho$ is not greater than the
diagonal of $P^m$, i.e. $\rho\leq\frac{\sqrt{S}}{2^m\sqrt{2}}$.
Therefore, using Lemma 3 and the mentioned inequalities for
$\delta(\zeta,\zeta_0)$ and $\rho$, we obtain
\begin{equation}\label{ocinka-4}
\Biggr\|\int\limits_{\partial
P^m_{\zeta}}\Psi(\zeta)\sigma\Biggr\|=\|I_5\|\leq 3nM\int\limits_{
\partial P^m_{\zeta}}\rho\,\|\delta(\zeta,\zeta_0)\|\,\|\sigma\|\leq
3nM\frac{\sqrt{S}}{2^m\sqrt{2}}\,\frac{S}{4^m}\,\varepsilon\,.
\end{equation}

It follows from the relations (\ref{ocinka-1}) and
(\ref{ocinka-4}) that\, $K\le c\,\varepsilon $, where the
constant\, $c$\, does not depend on\, $\varepsilon$. Passing to
the limit in the last inequality as $\varepsilon\rightarrow0$, we
obtain the equality $K=0$, and the theorem is proved.

\section{Main result}

Let us establish an analogue of Cauchy integral theorem for the
surface integral on the boundary $\partial\Omega_{\zeta}$ in the
case where the function
$\Psi:\overline{\Omega}_{\zeta}\rightarrow\mathbb{A}$ is
hyperholomorphic in a domain $\Omega_{\zeta}$ and continuous in
the closure $\overline{\Omega}_{\zeta}$ of this domain.

For such a function consider the modulus of continuity
$$\omega_{\overline{\Omega}_{\zeta}}(\Psi,\delta):=\sup\limits_{\zeta_1,\zeta_2\in
\overline{\Omega}_{\zeta},\|\zeta_1-\zeta_2\|\leq\delta}\|\Psi(\zeta_1)-\Psi(\zeta_2)\|.$$

The {\it two-dimensional upper Minkowski content} (see, e.g.,
\cite[p.~79]{Mattila}) is
$$
\mathcal{M}^*(\partial\Omega):=\limsup\limits_{\varepsilon\rightarrow0}\frac
{V(\partial\Omega^\varepsilon)}{2\varepsilon}\,,
$$
where $V(\partial\Omega^\varepsilon)$ denotes the volume of
$\partial\Omega^\varepsilon$.\vskip 2mm

{\bf Theorem 2.} \textit{Suppose that the boundary
$\partial\Omega$ of a simply connected domain
$\Omega\subset\mathbb{R}^{3}$ is a closed quadrable surface for
which $\mathcal{M}^*(\partial\Omega)<\infty$, and $\Omega$ has
Jordan measurable intersections with planes perpendicular to
coordinate axes. Suppose also that a function
$\Psi:\overline{\Omega}_{\zeta}\rightarrow\mathbb{A}$ is
hyperholomorphic in the domain $\Omega_{\zeta}$ and continuous in
the closure $\overline{\Omega}_{\zeta}$ of this domain. Then the
following equality holds:}
\begin{equation}\label{form-Koshi-po-pov}
\int\limits_{\partial \Omega_{\zeta}}\Psi(\zeta)\sigma=0.
\end{equation}

\textbf{\textit{Proof.}} Inasmuch as
$\mathcal{M}^*(\partial\Omega)<\infty$, there exists
$\varepsilon_0>0$ such that for all
$\varepsilon\in(0,\varepsilon_0)$ the following inequality holds:
\begin{equation}\label{teo2,ner-dlja-pover}
V(\partial\Omega^\varepsilon)\leq c\, \varepsilon,
\end{equation}
where the constant\, $c$\, does not depend on\, $\varepsilon$.

Let us take $\varepsilon<\varepsilon_0/\sqrt{3}$. Let us make a
partition of the space $\mathbb{R}^3$ onto cubes with an edge of
the length $\varepsilon$ by planes perpendicular to the coordinate
axes. Then we have the equality
\begin{equation}\label{form-Koshi-po-pov-1}
\int\limits_{\partial\Omega_{\zeta}}\Psi(\zeta)\,\sigma=\sum\limits_j
\int\limits_{\partial(\Omega_{\zeta}\cap
K^j_\zeta)}\Psi(\zeta)\,\sigma+\sum\limits_k \int\limits_{\partial
K^k_\zeta}\Psi(\zeta)\,\sigma,
\end{equation}
where the first sum is applied to the cubes $K^j$ for which
$\overline{K^j}\cap\partial\Omega\ne\varnothing$, and the second
sum is applied to the cubes $K^k$ for which
$\overline{K^k}\subset\Omega$. By Theorem 1, the second sum is
equal to zero.

To estimate an integral of the first sum we take a point
$\zeta_j\in \Omega_{\zeta}\cap K^j_\zeta$. Note that the diameter
of set $ \Omega\cap K^j$ does not exceed $ \varepsilon\sqrt {3}$.
Inasmuch as $\Omega$ has Jordan measurable intersections with
planes perpendicular to coordinate axes, the Lebesgue measure of
the boundaries of mentioned intersections is equal to $0$, and
consequently, the set $\partial(\Omega_{\zeta}\cap K^j_\zeta)$
consists of closed quadrable surfaces. Therefore, taking into
account the equality (\ref{int-sigma}) and using Lemma 3, we
obtain
$$\Biggr\|\int\limits_{\partial(\Omega_{\zeta}\cap
K^j_\zeta)}\Psi(\zeta)\sigma\Biggr\|=\Biggr\|\int\limits_{\partial(\Omega_{\zeta}\cap
K^j_\zeta)}(\Psi(\zeta)-\Psi(\zeta_j))\sigma\Biggr\|\leq$$
\begin{equation}\label{teo2_ocinka-1}
\leq 3nM\int\limits_{\partial(\Omega_{\zeta}\cap
K^j_\zeta)}\|\Psi(\zeta)-\Psi(\zeta_j)\| \|\sigma\| \leq
3nM\,\omega_{\overline{\Omega}_{\zeta}}(\Psi,\varepsilon\sqrt{3})\int\limits
_{\partial(\Omega_{\zeta}\cap K^j_\zeta)}\|\sigma\|.
\end{equation}

Thus, the following estimate is a result of the equality
(\ref{form-Koshi-po-pov-1}) and the inequality
(\ref{teo2_ocinka-1}):
\begin{equation}\label{teo2_ocinka}
\Biggr\|\int\limits_{\partial\Omega_\zeta}\Psi(\zeta)\,\sigma\Biggr\|
\leq3nM\omega_{\overline{\Omega}_{\zeta}}(\Psi,\varepsilon\sqrt{3})
\sum_j\int\limits_{\partial(\Omega_{\zeta}\cap
K^j_\zeta)}\|\sigma\|\leq$$
$$\leq3nM\omega_{\overline{\Omega}_{\zeta}}(\Psi,\varepsilon\sqrt{3})\bigg(\int\limits_
{\partial\Omega_{\zeta}}\|\sigma\|+6\sum_j\varepsilon^2\bigg).
\end{equation}

Inasmuch as $\bigcup_{j}K^j\subset\Omega^{\varepsilon\sqrt{3}}$,
taking into account the inequality (\ref{teo2,ner-dlja-pover}), we
obtain the estimation
$$\sum_j\varepsilon^3\leq V\Big(\partial\Omega^{\varepsilon\sqrt{3}}\Big)\le
c\varepsilon\sqrt{3}\,,$$
 from which it follows that
\begin{equation}\label{epsilon_kv}
\sum_j\varepsilon^2\leq c\sqrt{3}.
\end{equation}

Finally, the following inequality is as a result of the
estimations (\ref{teo2_ocinka}) and (\ref{epsilon_kv}):
\begin{equation}\label{epsilon_kv-1}
\Biggr\|\int\limits_{\partial\Omega_\zeta}\Psi(\zeta)\sigma\Biggr\|\leq
c_1\,\omega_{\overline{\Omega}_{\zeta}}(\Psi,\varepsilon\sqrt{3})
\end{equation}
where the constant\, $c_1$\, does not depend on\, $\varepsilon$.

To complete the proof, note that
$\omega_{\overline{\Omega}_\zeta}(\Psi,\varepsilon\sqrt{3})\to 0$
as $\varepsilon\rightarrow0$ due to the uniform continuity of the
function $\Psi$ on $\overline{\Omega}_{\zeta}$. \vskip 2mm

Theorem 2 generalizes Theorem 1 \cite{Pl-Shp3} that was proved in
a three-dimensional commutative algebra for functions which
generate solutions of the three-dimensional Laplace equation.

\section{Remarks}

Note that for a surface $\Sigma$ in $\mathbb{R}^3$, there exist
positive constants $c_1$ and $c_2$ such that
\begin{equation}\label{naslid-1}
c_1\varepsilon^3N_{\Sigma}(\varepsilon)\leq
V(\Sigma^\varepsilon)\leq
c_2\varepsilon^3N_{\Sigma}(\varepsilon)\,,
\end{equation}
where $N_{\Sigma}(\varepsilon)$ is the least number of
$\varepsilon$-balls needed to cover $\Sigma$ (see
\cite{Borodich}).

It is evidently follows from (\ref{naslid-1}) that the inequality
(\ref{teo2,ner-dlja-pover}) is equivalent to the inequality of the
form
\begin{equation}\label{naslid-2}
N_{\Sigma}(\varepsilon)\,\varepsilon^2\leq c,
\end{equation}
where the constant $c$ does not depend on $\varepsilon$.

Taking into account that a rectifiable surface $\Sigma$ is a
Lipschitz image of the square $G$ and the inequality of the form
(\ref{naslid-2}) is fulfilled for $G$, it is easy to prove the
inequality (\ref{naslid-2}) for $\Sigma$.

For a surface $\Sigma$ in $\mathbb{R}^3$ that has a finite
$2$-dimensional Hausdorff measure $\mathcal{H}^2(\Sigma)$, if
there exists a positive constant $c$ such that
\begin{equation}\label{left-reg}
c \varepsilon^2\leq\mathcal{H}^2\big(\Sigma\cap
B(x,\varepsilon)\big) \quad
\forall\,x\in\Sigma\,\,\,\forall\,\varepsilon\in(0;{\rm
diam}\,\Sigma]\,,
\end{equation}
where ${\rm diam}\,\Sigma$ is the diameter of $\Sigma$, and
$B(x,\varepsilon)$ denotes the open ball with center $x$ and
radius $\varepsilon$, then the inequalities
$P_\Sigma(\varepsilon)\varepsilon^2\leq
c_1\mathcal{H}^2(\Sigma)<\infty$ is fulfilled, where
$P_\Sigma(\varepsilon)$ is the greatest number of disjoint
$\varepsilon$-balls with centers in $\Sigma$ and the constant
$c_1$ does not depend on $\varepsilon$ (see
\cite[p.~309]{Abreu-Bory-Moreno}). Taking into account the
inequality $N_{\Sigma}(2\varepsilon)\leq P_\Sigma(\varepsilon)$
(see \cite[p.~78]{Mattila}), we obtain the inequality
(\ref{naslid-2}) for a surface $\Sigma$ satisfying the condition
(\ref{left-reg}). \vskip 2mm

The authors are grateful to Dr. M.~Tkachuk for very useful
discussions of geometric aspects of our researches.


\vspace{10mm}

Authors:\\[2mm]
{ \bf S.\,A.~Plaksa}, {Institute of Mathematics of the National
Academy of Sciences of Ukraine, Ukraine, 01601, Kiev,
Tereshchenkivska Str. 3,\\ Phone (office): (38044) 234 51 50,
E-mail: plaksa@imath.kiev.ua}\\[2mm]
{\bf V.\,S.~Shpakivskyi}, {Institute of Mathematics of the
National Academy of Sciences of Ukraine, Ukraine, 01601, Kiev,
Tereshchenkivska Str. 3,\\ Phone (office): (38044) 234 51 50,
E-mail: shpakivskyi@mail.ru}

\end{document}